\documentclass[11pt,a4paper]{article}
\usepackage{mathrsfs}
\usepackage{amssymb}
\usepackage{amsfonts}
\usepackage{color,xcolor}
\usepackage{graphicx}
\usepackage{manfnt}
\RequirePackage[colorlinks,citecolor=blue,urlcolor=blue]{hyperref}
\newtheorem{theorem}{Theorem}[section]

\newtheorem{lemma}{Lemma}[section]

\newtheorem{definition}{Definition}[section]

\def\<{{\langle}}\def\>{{\rangle}}
\def\[{{\Big[}}\def\]{{\Big]}}\def\({{\Big(}}\def\){{\Big)}}

\def\={&\!\!=\!\!&}
\def\cA{{\mathcal A}}

\def\cS{{\mathcal S}}

\def\mN{{\mathbb N}}

\def\mZ{{\mathbb Z}}

\def\geq{\geqslant}\def\leq{\leqslant}

\def\div{\mathord{{\rm div}}}

\begin{document}
\title{\bf Strong solution for stochastic transport equations with irregular drift: existence and non-existence\footnote{This work was partly supported by the NSFC grants 11501577, 11531006, 11771449, 11771123 and PAPD of Jiangsu Higher Education Institutions.}}

\author{Jinlong Wei$^a$, Jinqiao Duan$^b$, Hongjun Gao$^c$ and Guangying Lv$^d$}
\date{{$^a$ School of Statistics and Mathematics, Zhongnan University of}\\
{Economics and Law, Wuhan, Hubei 430073, China}\\
{weijinlong@zuel.edu.cn}
\\ {$^b$ Department of Applied Mathematics}\\
{Illinois Institute of Technology, Chicago, IL 60616, USA}\\
{duan@iit.edu} \\  {$^c$ Institute of Mathematics, School of Mathematical Science} \\ {Nanjing Normal University, Nanjing 210023, China}\\ {gaohj@njnu.edu.cn}
\\ {$^d$ Institute of Contemporary Mathematics, Henan University}\\
{Kaifeng, Henan 475001, China} \\
{gylvmaths@henu.edu.cn}}

 \maketitle
\noindent{\hrulefill}

\vskip1mm\noindent{\bf Abstract} We prove some existence, uniqueness and non-existence results of stochastic strong solutions for a class of stochastic transport equations with a $q$-integrable (in time), bounded and $\alpha$-H\"{o}lder continuous (in space) drift coefficient. More precisely, we show that for a Sobolev differentiable initial condition, there exists  a unique stochastic strong solution when $\alpha>2/q$, while   for $\alpha+1<2/q$ with spatial dimension higher than one, we can choose  proper initial data and drift coefficients so that there is no stochastic strong solutions.

 \vskip1.2mm\noindent
{\bf MSC (2010):} 60H15 (35A01 35L02)

\vskip1.2mm\noindent
{\bf Keywords:} Stochastic transport equations; Stochastic strong solution; Existence; Nonexistence; Stochastic PDEs

 \vskip0mm\noindent{\hrulefill}
\section{Introduction}\label{sec1}\setcounter{equation}{0}
Consider the Cauchy problem
\begin{eqnarray}\label{1.1}
\left\{
  \begin{array}{ll}
du(t,x)+b(t,x)\cdot\nabla u(t,x)dt
+\sum_{i=1}^d\partial_{x_i}u(t,x)\circ dB_i(t)=0, \ (\omega,t,x)\in \Omega\times (0,T)\times {\mathbb R}^d, \\
u(0,x)=u_0(x), \  x\in{\mathbb R}^d,
  \end{array}
\right.
\end{eqnarray}
where $T>0$ is a given real number, $B(t)=(B_1(t), B_2(t), _{\cdots}, B_d(t))$ is a $d$-dimensional standard Brownian motion defined on a stochastic basis ($\Omega, \mathcal{F},{\mathbb P},(\mathcal{F}_{t})_{t\geq 0}$), and the stochastic
integration with a notation $\circ$ is interpreted in Stratonovich
sense. The drift coefficient $b: [0,T]\times{\mathbb R}^d\rightarrow{\mathbb R}^d$ and the initial value $ u_0:
{\mathbb R}^d\rightarrow{\mathbb R}$ are measurable functions    $L^1(0,T;L^1_{loc}({\mathbb R}^d;{\mathbb R}^d))$ and $L^1(0,T;L^1_{loc}({\mathbb R}^d))$, respectively. We will prove some existence and non-existence results on stochastic strong solutions. Here a stochastic strong solution is defined as  follows \cite{CO,FF2,FGP1}.
\begin{definition} \label{def1.1} Let $p\in [1,\infty]$ and $\div b\in L^1(0,T;L^{p^\prime}_{loc}({\mathbb R}^d))$ with $1/p+1/p^\prime=1$.   A random field $u$ in $ L^\infty(\Omega\times[0,T];L^p({\mathbb R}^d))$
  is called a stochastic weak solution of (\ref{1.1}) if for every
$\varphi\in\mathcal{C}_0^\infty({\mathbb R}^d)$, $\int_{{\mathbb R}^d}\varphi(x)u(t,x)dx$
has a continuous modification which is an $\mathcal{F}_t$-semimartingale,  and
the following identity holds,
\begin{eqnarray}\label{1.2}
\int_{{\mathbb R}^d}\varphi(x)u(t,x)dx&=&\int_{{\mathbb R}^d}\varphi(x)u_0(x)dx+
\int^t_0\int_{{\mathbb R}^d}\div(b(s,x)\varphi(x))u(s,x)dxds\cr\cr&&
+\sum_{i=1}^d\int^t_0\circ dB_i(s)\int_{{\mathbb R}^d}\partial_{x_i}\varphi(x)u(s,x)dx,  \quad
{\mathbb P}-a.s. \;\; t\in [0,T].
\end{eqnarray}
Moreover, if the following additional estimates hold
\begin{eqnarray}\label{1.3}
\left\{
  \begin{array}{ll}
{\mathbb E}\|\nabla u\|_{L^\infty(0,T;L^p_{loc}({\mathbb R}^d))}^p<\infty, \quad when \ p<\infty, \\
{\mathbb E}\|\nabla u\|_{L^\infty(0,T;L^\infty_{loc}({\mathbb R}^d))}^r<\infty, \quad \forall \ r\in [1,\infty),
  \end{array}
\right.
\end{eqnarray}
then $u$ is called  a stochastic strong solution.
\end{definition}

There is a great recent interest in studying the existence and uniqueness of stochastic weak solutions in (\ref{1.1}). When $p=\infty$,  this is investigated in \cite{FGP1} if the drift coefficient $b$ is H\"{o}lder continuous in spatial variable and bounded in temporal variable, and  in  \cite{NO} if $b$ is only integrable in spatial and temporal variables. Some other extensions for $p=\infty$ have also been established  \cite{AF,FL,MNP,Zha}. For $p<\infty$, the existence and uniqueness of stochastic weak solution,  in the case of a Sobolev differentiable drift coefficient $b$, is shown by  Catuogno and Olivera \cite{CO}.

The well-posedness of stochastic weak solutions   in $\cap_{r\geq 1}W^{1,r}_{loc}({\mathbb R}^d)$  have also been examined. When the drift coefficient $b$ is only integrable in space and time, which satisfies a Ladyzhenskaya-Prodi-Serrin condition, i.e.
\begin{eqnarray}\label{1.4}
b\in L^q(0,T;L^p({\mathbb R}^d;{\mathbb R}^d)), \ q,p\in [2,\infty), \ \frac{2}{q}+\frac{d}{p}\leq 1, \ or \ p=\infty, \ q=2,
\end{eqnarray}
for $u_0\in \cap_{r\geq 1}W^{1,r}({\mathbb R}^d)$, the existence and uniqueness of such stochastic weak solutions   is shown in \cite{BFG,FF2}. However, if the initial condition  is merely assumed to be Sobolev differentiable   in $W^{1,p}({\mathbb R}^d)$ for a fixed $p\in[1,\infty]$, the existence and uniqueness of    stochastic weak solutions  in $W^{1,p}_{loc}({\mathbb R}^d)$   is still open.  But when the Ladyzhenskaya-Prodi-Serrin condition is replaced by $b\in L^\infty(0,T;\mathcal{C}_b^\beta({\mathbb R}^d;{\mathbb R}^d))$ ($\beta\in (0,1)$),  Flandoli, Gubinelli and Priola  \cite{FGP3}  give a positive answer to this open problem.  Hence two related interesting questions  may be asked:

$\bullet$ Does there exist a unique stochastic weak solution   in $W^{1,p}_{loc}({\mathbb R}^d)$, almost surely,  for drift coefficient $b$  integrable in time, bounded and $\alpha$-H\"{o}lder continuous in space?

$\bullet$ Is there an integrable drift coefficient $b$ which does not satisfy a Ladyzhenskaya-Prodi-Serrin type condition,  together with a $W^{1,p}({\mathbb R}^d)$ initial condition, such that there is no-existence of  a $W^{1,p}_{loc}({\mathbb R}^d)$ solution?

The purpose  of this paper is to answer both questions. To this end,  we introduce a concept of stochastic strong solution (Definition \ref{def1.1}) and  give two positive answers for both questions.

The  answer to the first question is described by the following theorem.

\begin{theorem} \label{the1.1} \textbf{(Existence and uniqueness)} Assume that $q\in [1,\infty]$, $\alpha\in (0,1)$ and $p\in [1,\infty]$ such that
\begin{eqnarray}\label{1.5}
b\in L^q(0,T;\mathcal{C}^\alpha_b({\mathbb R}^d;{\mathbb R}^d)), \ \ div b\in L^1(0,T;L^\infty({\mathbb R}^d)), \ \ u_0\in W^{1,p}({\mathbb R}^d).
\end{eqnarray}
If $2/q<\alpha$,  then there exists a unique stochastic strong solution to the Cauchy problem (\ref{1.1}).  Moreover, the unique
stochastic strong solution can be represented by $u(t,x)=u_0(X^{-1}(t,x))$, with $\{X(t,x)\}$ being the unique
strong solution of the associated stochastic differential equation  (\ref{2.1}) with $s=0$.
 \end{theorem}

\vskip2mm\noindent
\textbf{Remark 1.1.} (i) When $q=\infty$ and (\ref{1.5}) holds, there exists a unique $W^{1,p}_{loc}$-solution $u$ of the Cauchy problem (\ref{1.1}) by Flandoli, Gubinelli and Priola \cite[Theorems 5 and 16]{FGP1},  i.e. $u$ is a stochastic weak solution and  $u(\omega,t)\in W^{1,p}_{loc}$ almost surely for every $t\in[0,T]$. Here we generalize the result  of Flandoli, Gubinelli and Priola \cite{FGP1} to the case of $q\in [1,\infty]$. Moreover, for $q\in [1,\infty]$, we may get an analogue result of \cite[Theorem 20]{FGP1} for stochastic weak solutions in the case of $p=\infty$.

(ii) When  the noise vanishes,  we \cite{WDGL} have constructed a counterexample to illustrate the non-existence of $W^{1,p}_{loc}$ solutions for the corresponding deterministic equation. So this theorem shows that the noise can regularize the solutions.

(iii)  Several other works have previously explained that a noise can regularize the solutions. For example, noise can prevent infinite stretching of the passive field in a stochastic vector advection equation \cite{FMN}, or can prevent collapse of Vlasov-Poisson point charges \cite{DFV}.

For the second question, the answer is in the following theorem.

\begin{theorem} \label{the1.2} \textbf{(Non-existence)} Assume that  $q\in [1,\infty]$, $\alpha\in (0,1)$ and $p\in [1,\infty]$ such that (\ref{1.5})  holds.  If $2/q>\alpha+1$ and $d\geq 2$,  and for every $p\in [1,\infty]$, one can choose a proper initial condition  $u_0$ and a drift coefficient $b$ such that there exists a unique stochastic weak solution to the Cauchy problem (\ref{1.1}) which can be represented by $u(t,x)=u_0(X^{-1}(t,x))$. However when $p\in [1,\infty)$, $u$ does not lie in
$L^p(\Omega;L^\infty(0,T;W^{1,p}_{loc}({\mathbb R}^d)))$ and when $p=\infty$, $u$ does not lie in
$L^r(\Omega;L^\infty(0,T;W^{1,\infty}_{loc}({\mathbb R}^d)))$, for every $r\geq 1$.
 \end{theorem}

\vskip2mm\noindent
\textbf{Remark 1.2.} (i)  When the noise vanishes,   (\ref{1.1}) reduces to
\begin{eqnarray}\label{1.6}
\left\{
  \begin{array}{ll}
\partial_tu(t,x)+b(t,x)\cdot\nabla u(t,x)=0, \ (t,x)\in (0,T)\times {\mathbb R}^d, \\
u(t,x)|_{t=0}=u_0(x), \  x\in{\mathbb R}^d.
  \end{array}
\right.
\end{eqnarray}
For a particular  drift coefficient $b$ (given in  Section \ref{sec4}), there exists a unique weak solution to (\ref{1.6}). Moreover, the unique weak
solution can be represented by $u(t,x)=u_0(X^{-1}(t,x))$. Here
$X(t,x)$ is the unique solution of the ordinary differential equation (ODE) $\dot{X}=b(t,X)$. However, $u(t,x)=u_0(X^{-1}(t,x))$ does not lie in
$L^\infty(0,T;W^{1,p}_{loc}({\mathbb R}^d)).$ Here $u$ is said to be a weak
solution of (\ref{1.6}), if it lies in $L^\infty(0,T;L^p({\mathbb R}^d))$
and for every $\varphi\in\mathcal{C}_0^\infty({\mathbb R}^d)$,
\begin{eqnarray*}
\int_{{\mathbb R}^d}\varphi(x)u(t,x)dx=\int_{{\mathbb R}^d}\varphi(x)u_0(x)dx+
\int^t_0\int_{{\mathbb R}^d}\div(b(s,x)\varphi(x))u(s,x)dxds,
\end{eqnarray*}
for every $t\in [0,T]$.    If additionally $|\nabla u|\in L^\infty(0,T;L^p_{loc}({\mathbb R}^d))$, we call $u$ a strong solution of (\ref{1.6}). In this sense, the noise has no regularizing effect.

(ii)  Assume that $p\in [1,\infty]$ and $k\in L^p({\mathbb R}^d)$, such that  for every $\tau>0$,
\begin{eqnarray}\label{1.7}
\Big[\int_{{\mathbb R}^d}|k(\tau x)|^pdx\Big]^{\frac{1}{p}}=
\tau^{-\frac{d}{p}}\Big[\int_{{\mathbb R}^d}|k(x)|^pdx\Big]^{\frac{1}{p}}.
\end{eqnarray}
Defining $degree(k)=-d/p$ and noticing that for second order parabolic equations we can  `trade' space-regularity against time-regularity at a cost of one time derivative
for two space derivatives, we obtain that $degree(b)=-2/q-d/p$ if $b\in L^q(0,T;L^p({\mathbb R}^d;{\mathbb R}^d))$.  From this viewpoint, the  Ladyzhenskaya-Prodi-Serrin condition can be read as $degree(b)\leq -1$. Similarly, if $b\in L^q(0,T;\mathcal{C}^\alpha_b({\mathbb R}^d;{\mathbb R}^d))$, we derive an analogue of Ladyzhenskaya-Prodi-Serrin condition:
\begin{eqnarray}\label{1.8}
b\in L^q(0,T;\mathcal{C}^\alpha_b({\mathbb R}^d;{\mathbb R}^d)),  \ \frac{2}{q}\leq \alpha+ 1.
\end{eqnarray}
We call (\ref{1.8}) a  `Ladyzhenskaya-Prodi-Serrin type condition'. Thus, we give a positive answer for the second question.

(iii) When (\ref{1.8}) holds,  we may get an analogous  result as in \cite{BFG,FF2}. Precisely speaking, if $u_0\in \cap_{r\geq 1}W^{1,r}({\mathbb R}^d)$, there exists a unique stochastic weak solution $u$ to  (\ref{1.1}). Besides,  for every $t\in [0,T]$,
\begin{eqnarray*}
{\mathbb P}\Big(\ |\nabla u(t)|\in \cap_{r\geq 1}L^r_{loc}({\mathbb R}^d)\Big)=1.
\end{eqnarray*}

\medskip

The rest of this paper is arranged as follows. In Sections 2, we present some new results on existence and uniqueness of stochastic flow of diffeomorphisms on stochastic differential equations (SDEs).  In Sections 3 and 4, the
proof of Theorems \ref{the1.1} and \ref{the1.2} are given, respectively. Finally in the  Appendix,  we prove a useful lemma  needed for proving Theorem \ref{the1.1}.

\vskip2mm\noindent
\textbf{Notations}  The letter $C$ will mean a positive constant, whose values may change in different places. For a parameter or a function $\varrho$, $C(\varrho)$ means the constant is only dependent on $\varrho$. $\mN$ is the set of natural numbers and $\mZ$ denotes the set of integral numbers. For every $R>0$, $B_R:=\{x\in{\mathbb R}^d:|x|<R\}$. Almost surely can be abbreviated to $a.s.$.

\section{An SDE driven by an $L^q(0,T;\mathcal{C}_b^\alpha({\mathbb R}^d;{\mathbb R}^d))$ drift}\label{sec2}\setcounter{equation}{0}

Given $s\in[0,T]$ and $x\in {\mathbb R}^d$, consider the  following stochastic differential equation (SDE) in ${\mathbb R}^d$:
\begin{eqnarray}\label{2.1}
dX(s,t)=b(t,X(s,t))dt+dB(t), \ \ t\in(s,T], \ \ X(s,t)|_{t=s}=x.
\end{eqnarray}

\begin{definition} [\cite{Kun2}, P114] \label{def2.1} A stochastic homeomorphism flow (respect. of class $\mathcal{C}^{1,\beta}$ with $\beta\in (0,1)$) on
$(\Omega, \mathcal{F},{\mathbb P}, (\mathcal{F}_t)_{0\leq t\leq T})$ associated to (\ref{2.1}) is a
map $(s,t,x,\omega) \rightarrow X(s,t,x)(\omega)$, defined for
$0\leq s \leq t \leq T, \ x\in {\mathbb R}^d, \ \omega \in \Omega$ with
values in ${\mathbb R}^d$, such that

(i) given every  $s\in [0,T],\  x \in {\mathbb R}^d$, the process
$\{X(s,\cdot,x)\}= \{X(s,t,x), \ t\in [s,T]\}$ is a continuous
$\{\mathcal{F}_{s,t}\}_{s\leq t\leq T}$-adapted solution of (\ref{2.1});

(ii) ${\mathbb P}-a.s.$, for all  $0\leq s \leq t \leq T$, the functions
$X(s,t,x), \ X^{-1}(s,t,x)$ are continuous in $(s,t,x)$;

(iii) ${\mathbb P}-a.s., \   X(s,t,x)=X(r,t,X(s,r,x))$  for all
$0\leq s\leq r \leq t \leq T$, $x\in {\mathbb R}^d$ and $X(s,s,x)=x$.
\end{definition}

We now present an important   result.

\begin{lemma} \label{lem2.1}  Assume that $q\in [1,\infty]$, $\alpha \in (0,1)$ and  $2/q<\alpha$ such
that $b\in L^q(0,T; \mathcal{C}^\alpha_b({\mathbb R}^d;{\mathbb R}^d))$. Then for
every $s\in [0,T ]$ and $ x\in {\mathbb R}^d$, the stochastic differential equation
(\ref{2.1}) has a unique continuous adapted solution $\{X(s,t,x), \
t\in [s,T ], \ \omega \in\Omega\}$, which forms a $\mathcal{C}^{1,\alpha^\prime}$ ($\alpha^\prime<\alpha-2/q$) stochastic flow $X(s,t)$ of diffeomorphisms.
\end{lemma}

Before giving the proof, we need a lemma.
\begin{lemma} \label{lem2.2} Let $q\in[1,\infty]$, $\alpha\in (0,1)$ and $b,f\in L^q(0,T;\mathcal{C}^\alpha_b({\mathbb R}^d;{\mathbb R}^d))$ with $2/q<\alpha$.  For a real number $\lambda>0$, consider the following backward heat equation
\begin{eqnarray}\label{2.2}
\left\{
\begin{array}{ll}
\partial_{t}U(t,x) +\frac{1}{2}\Delta U(t,x)+b(t,x)\cdot \nabla U(t,x)=
\lambda U(t,x)+f(t,x), \ (t,x)\in (0,T)\times {\mathbb R}^d, \\
U(T,x)=0, \  x\in{\mathbb R}^d.
  \end{array}
\right.
\end{eqnarray}
We have the following assertions:

(i) there is a unique $U\in L^q(0,T;\mathcal{C}^{2,\alpha}_b({\mathbb R}^d;{\mathbb R}^d))\cap W^{1,q}(0,T;\mathcal{C}^\alpha_b({\mathbb R}^d;{\mathbb R}^d))$ solving the Cauchy problem (\ref{2.2});

(ii) $U\in L^\infty(0,T;\mathcal{C}^{2,\alpha-2/q}_b({\mathbb R}^d;{\mathbb R}^d))$ and there is a constant $C>0$, such that
\begin{eqnarray}\label{2.3}
\|U\|_{L^\infty(0,T;\mathcal{C}^{2,\alpha-2/q}_b({\mathbb R}^d;{\mathbb R}^d))}\leq C \|f\|_{L^q(0,T;\mathcal{C}^\alpha_b({\mathbb R}^d;{\mathbb R}^d))};
\end{eqnarray}

(iii) as $\lambda$ tends to 0,
\begin{eqnarray}\label{2.4}
\|\nabla U\|_{L^\infty([0,T]\times{\mathbb R}^d)}\longrightarrow 0.
\end{eqnarray}
\end{lemma}

\vskip2mm\noindent
\textbf{Proof.} Assertion (i) can be seen in \cite[Theorem 2.1]{Kry2}. It remains to check (ii) and (iii). Here we only examine the simple case $b=0$.  See \cite{KP} for details for general $b$ by using a continuity method.  First, let us calculate (ii). When $b=0$, $U$ has the following obvious representation
\begin{eqnarray}\label{2.5}
U(t,x)=\int_0^{T-t}e^{-\lambda r}P_{r}f(t+r,\cdot)(x)dr,
\end{eqnarray}
where $P_r$ is defined by
\begin{eqnarray}\label{2.6}
P_r\varphi(x)=\frac{1}{(2\pi r)^{d/2}}\int_{{\mathbb R}^d}e^{-\frac{|x-y|^2}{2r}}\varphi(y)dy:=(K(r,\cdot)\ast\varphi)(x), \quad \varphi\in L^\infty({\mathbb R}^d).
\end{eqnarray}
By (\ref{2.5}) and the assumption   $2/q<\alpha$,  we then know that  $U\in L^\infty(0,T;\mathcal{C}^2_b({\mathbb R}^d;{\mathbb R}^d))$ and there is a positive constant $C>0$ such that
\begin{eqnarray}\label{2.7}
\|U\|_{L^\infty(0,T;\mathcal{C}^2_b({\mathbb R}^d;{\mathbb R}^d))}\leq C \|f\|_{L^q(0,T;\mathcal{C}^\alpha_b({\mathbb R}^d;{\mathbb R}^d))}.
\end{eqnarray}

It needs to show that $\partial_{x_i}^2U\in L^\infty(0,T;\mathcal{C}^{\alpha-2/q}_b({\mathbb R}^d;{\mathbb R}^d))$ for every $1\leq i\leq d$,    and that  (\ref{2.3}) holds. For every $x,y\in{\mathbb R}^d$ and $1\leq i\leq d$,
\begin{eqnarray*}
&&\partial_{x_i}^2U(t,x)-\partial_{y_i}^2U(t,y)\cr\cr
&=&\int_0^{T-t}dr\int_{|x-z|\leq 2|x-y|}e^{-\lambda r}\partial_{x_i}^2K(r,x-z)[f(t+r,z)-f(t+r,x)]dz
\cr\cr&&-\int_0^{T-t}e^{-\lambda r}dr\int_{|x-z|\leq 2|x-y|}\partial_{y_i}^2K(r,y-z)[f(t+r,z)-f(t+r,y)]dz
\cr\cr&&+\int_0^{T-t}e^{-\lambda r}dr\int_{|x-z|> 2|x-y|}\partial_{y_i}^2K(r,y-z)[f(t+r,y)-f(t+r,x)]dz
\cr\cr&&+\int_0^{T-t}e^{-\lambda r}dr\int_{|x-z|> 2|x-y|}[\partial_{x_i}^2K(r,x-z)-\partial_{y_i}^2K(r,y-z)][f(t+r,z)-f(t+r,x)]dz
\cr\cr&&\cr\cr&=&:I_1(t)+I_2(t)+I_3(t)+I_4(t).
\end{eqnarray*}

Let us calculate $I_1, I_2, I_3, I_4$,  respectively. To start with, we manipulate the term $I_1$.
\begin{eqnarray}\label{2.8}
|I_1(t)|&\leq&\int_0^{T-t}\int_{|x-z|\leq 2|x-y|}|\partial_{x_i}^2K(r,x-z)||f(t+r,z)-f(t+r,x)|dzdr
\cr\cr&\leq& \int_0^{T-t}\int_{|x-z|\leq 2|x-y|} r^{-\frac{d+2}{2}}e^{-\frac{|x-z|^2}{2r}}|x-z|^\alpha [f]_{\alpha}(t+r)dzdr.
\end{eqnarray}
By utilizing the H\"{o}lder inequality, from (\ref{2.8}), it yields that
\begin{eqnarray}\label{2.9}
|I_1(t)|&\leq&C\|f\|_{L^q(0,T;\mathcal{C}_b^\alpha({\mathbb R}^d))} \int_{|x-z|\leq 2|x-y|} \Big|\int_0^{T-t}r^{-\frac{(d+2)q^\prime}{2}}e^{-\frac{q^\prime|x-z|^2}{2r}} dr\Big|^{\frac{1}{q^\prime}} |x-z|^\alpha dz
\cr\cr&\leq& C\|f\|_{L^q(0,T;\mathcal{C}_b^\alpha({\mathbb R}^d))} \int_{|x-z|\leq 2|x-y|} \Big|\int_0^\infty r^{-\frac{(d+2)q^\prime}{2}}e^{-\frac{q^\prime|x-z|^2}{2r}} dr\Big|^{\frac{1}{q^\prime}} |x-z|^\alpha dz
\cr\cr&=& C\|f\|_{L^q(0,T;\mathcal{C}_b^\alpha({\mathbb R}^d))} \int_{|x-z|\leq 2|x-y|} \Big|\int_0^\infty r^{\frac{(d+2)q^\prime}{2}-2}e^{-\frac{q^\prime r}{2}} dr\Big|^{\frac{1}{q^\prime}} |x-z|^{d-2+\alpha+\frac{2}{q^\prime}} dz
\cr\cr&=& C\|f\|_{L^q(0,T;\mathcal{C}_b^\alpha({\mathbb R}^d))} |x-y|^{\alpha+\frac{2}{q^\prime}-2}
\cr\cr&=& C\|f\|_{L^q(0,T;\mathcal{C}_b^\alpha({\mathbb R}^d))} |x-y|^{\alpha-\frac{2}{q}},
\end{eqnarray}
where in the last identity we have used $q^\prime=q/(q-1)$.

Similarly, one gets that
\begin{eqnarray}\label{2.10}
|I_2(t)|\leq C\|f\|_{L^q(0,T;\mathcal{C}_b^\alpha({\mathbb R}^d))} |x-y|^{\alpha-\frac{2}{q}}.
\end{eqnarray}

For $I_3$, we employ Gauss-Green's formula primarily to gain
\begin{eqnarray}\label{2.11}
I_3(t)=\int_0^{T-t}e^{-\lambda r}
dr\int_{|y-z|=2|x-y|}\partial_{y_i}K(r,y-z)n_i[f(t+r,y)-f(t+r,x)]dS.
\end{eqnarray}
From (\ref{2.11}), owing to the H\"{o}lder inequality, we arrive at
\begin{eqnarray}\label{2.12}
&&|I_3(t)|\cr\cr&\leq& C\|f\|_{L^q(0,T;\mathcal{C}_b^\alpha({\mathbb R}^d))}
|x-y|^{\alpha}\int_{|x-z|=2|x-y|} \Big(\int_0^\infty r^{-\frac{q^\prime (d+1)}{2}}e^{-\frac{q^\prime|y-z|^2}{2r}}dr\Big)^{\frac{1}{q^\prime}}dS
\cr\cr&\leq&C\|f\|_{L^q(0,T;\mathcal{C}_b^\alpha({\mathbb R}^d))} |x-y|^{\alpha}\int_{|x-z|=2|x-y|} |y-z|^{-d-1+\frac{2}{q^\prime}}dS\Big(\int_0^\infty r^{\frac{q^\prime (d+1)}{2}-2}e^{-\frac{(d+1)q^\prime r}{2}}dr\Big)^{\frac{1}{q^\prime}}
\cr\cr\cr\cr&\leq&C\|f\|_{L^q(0,T;\mathcal{C}_b^\alpha({\mathbb R}^d))}|x-y|^{\alpha-\frac{2}{q}}.
\end{eqnarray}

To calculate $I_4$, we also use the H\"{o}lder inequality, and then acquire
\begin{eqnarray*}
&&|I_4(t)|\cr\cr &\leq& C\|f\|_{L^q(0,T;\mathcal{C}_b^\alpha({\mathbb R}^d))}
\int_{|x-z|> 2|x-y|}|x-z|^\alpha\Big(\int_0^{T-t}|\partial^2_{x_i}
K(r,x-z)-\partial_{y_i}^2K(r,y-z)|^{q^\prime}dr
\Big)^{\frac{1}{q^\prime}}dz.
\end{eqnarray*}
Notice that $|x-z|>2|x-y|$. So for every $\xi\in [x,y]$,
$$
\frac{1}{2}|x-z| \leq |\xi-z|\leq 2|x-z|.
$$
By virtue of mean value inequality, we have
\begin{eqnarray}\label{2.13}
&&|I_4(t)|\cr\cr &\leq& C\|f\|_{L^q(0,T;\mathcal{C}_b^\alpha({\mathbb R}^d))}|x-y|\int_{|x-z|> 2|x-y|}|x-z|^\alpha \Big(\int_0^{T-t} r^{-\frac{(d+3)q^\prime}{2}}
e^{-\frac{q^\prime|x-z|^2}{8r}}dr\Big)^{\frac{1}{q^\prime}}dz\cr\cr&\leq&
C\|f\|_{L^q(0,T;\mathcal{C}_b^\alpha({\mathbb R}^d))} |x-y|\int_{|x-z|> 2|x-y|}|x-z|^{\alpha-d-3+\frac{2}{q^\prime}} \Big(\int_0^\infty
r^{\frac{(d+3)q^\prime}{2}-2}e^{-\frac{q^\prime r}{8}}dr\Big)^{\frac{1}{q^\prime}}dz \cr\cr&\leq&C\|f\|_{L^q(0,T;\mathcal{C}_b^\alpha({\mathbb R}^d))}|x-y|^{\alpha-\frac{2}{q}}.
\end{eqnarray}

Combining (\ref{2.7}), (\ref{2.9})-(\ref{2.10}), (\ref{2.12})-(\ref{2.13}), then (\ref{2.3}) holds true. So we finish the proof for assertion (ii).

To check assertion (iii), by using the explicit formula (\ref{2.5}) for every $1\leq i\leq d$, then
\begin{eqnarray*}
|\partial_{x_i}U(t,x)|&=&\Big|\int_0^{T-t}dr\int_{{\mathbb R}^d}e^{-\lambda r}\partial_{x_i}K(r,x-z)[f(t+r,z)-f(t+r,x)]dz\Big|
\cr\cr&\leq&\int_0^{T-t}[f]_\alpha(t+r)e^{-\lambda r} dr \int_{{\mathbb R}^d}
r^{-\frac{d+1}{2}}e^{-\frac{|z|^2}{2r}}|z|^\alpha dz
cr\cr\cr&\leq&C\int_0^{T-t}[f]_\alpha(t+r)r^{-\frac{1-\alpha}{2}}e^{-\lambda r} dr\cr\cr&\leq&C\|f\|_{L^q(0,T;\mathcal{C}_b^\alpha({\mathbb R}^d))}
\Big(\int_0^{T-t}r^{-\frac{(1-\alpha)q}{(q-2)}} dr\Big)^{\frac{q-2}{2q}}\Big(\int_0^{T-t}e^{-2\lambda r} dr\Big)^{\frac{1}{2}}
\cr\cr&\leq&C\|f\|_{L^q(0,T;\mathcal{C}_b^\alpha({\mathbb R}^d))}\lambda^{-\frac{1}{2}},
\end{eqnarray*}
where in the last inequality we have used the fact $\alpha>2/q$. So this completes the proof. $\Box$

\vskip2mm\noindent
\textbf{Remark 2.1.} The idea to introduce Lemma \ref{lem2.2} is enlightened by \cite[Theorem 2]{FGP1}, which will serve us well later in proving Lemma \ref{lem2.1}. The use of Lemma \ref{lem2.2} is not limited to that extent. As we all know, for a second parabolic equation (such as (\ref{2.2})) with $b,f$   bounded (from Remark 1.2 (ii) now $degree(f)=degree(b)=0$), in general $U$ does  not lie in $W^{2,\infty}$. But when
$b,f$ are bounded and $\varsigma$-H\"{o}lder continuous in space ($degree(f)=degree(b)>0$), we have that  $U\in L^\infty(0,T;C_{b}^{2,\varsigma})$. The advantage of the present lemma is that it helps  derive a $W^{2,\infty}$ estimate for solutions in critical spaces in the sense of $degree(f)=degree(b)=0$.


\vskip2mm\noindent
\textbf{Proof of Lemma \ref{lem2.1}.} We only recall the idea of the proof; see \cite[Theorem 5]{FGP1} for details. Let $U$ be the unique solution of (\ref{2.2}) with the nonhomogeneous term $f$ is replaced by $-b$ and we define $\gamma(t,x)=x+U(t,x)$. For $\lambda$ sufficiently large, $\gamma(t)$ forms a non-singular diffeomorphism of class $\mathcal{C}^2$ uniformly in $t\in [0,T]$. Besides, for every $t\in [0,T]$, the inverse of $\gamma(t)$ (denoted by $\gamma^{-1}(t)$ has bounded first and second spatial derivatives, uniformly in $t\in [0,T]$. Consider the following SDE (see \cite{FF2,FF3,FGP1}):
\begin{eqnarray}\label{2.14}
d Y(t)=\lambda U(t,\gamma^{-1}(t,Y(t)))dt+ [I+\nabla
U(t,\gamma^{-1}(t,Y(t)))] dB(t),  t\in(s,T], \ Y(t)|_{t=s}=y.
\end{eqnarray}
Then it is equivalent to (\ref{2.1}). Since now (\ref{2.3}) holds, the classical results (see \cite[Chap. 2]{Kun1})  imply  the existence and uniqueness of a $\mathcal{C}^{1,\alpha^\prime}$ ($\alpha^\prime<\alpha-2/q$) stochastic flow of diffeomorphisms of (\ref{2.14}). On the other hand, the relationship between (\ref{2.1}) and (\ref{2.14}) is given by $X(t)=\gamma^{-1}(t,Y(t))$.  This finishes the proof. $\Box$

\vskip2mm\noindent
\textbf{Remark 2.2.}  When $q=\infty$, Lemma \ref{lem2.1} degenerates into \cite[Theorem 5]{FGP1}, but when $q<\infty$ this result is new.  Moreover, when $2/q<\alpha+1$ we also gain the existence and uniqueness of $\mathcal{C}^{\alpha^\prime}$ ($\alpha^\prime<1$) stochastic flow of homeomorphism. For more details in the case of $q=\infty$, one can refers to \cite{Att,FF1,FF3,Fl,FGP2,KR} and the references cited  therein. For more details for SDEs, one consults to \cite{Duan}.

\section{Proof of Theorem \ref{the1.1}}\label{sec3}
\setcounter{equation}{0}
\vskip2mm\noindent
\textbf{Proof.} First, we check the uniqueness and observing that the equation is linear, it suffices to prove that $u\equiv 0$ a.s. if the initial data vanishes. Let $\varrho_\varepsilon$ be a regularizing kernel i.e.
\begin{eqnarray*}
\varrho_\varepsilon =\frac{1}{\varepsilon^d} \varrho(\frac{\cdot}{\varepsilon}) \ \ with \ \ 0\leq \varrho \in \mathcal{C}^\infty_0({\mathbb R}^d) , \ \ support(\varrho)\subset B_1, \ \ \int_{{\mathbb R}^d}\varrho(x)dx=1.
\end{eqnarray*}
We set $u_\varepsilon=u\ast \varrho_\varepsilon$. Then $u_\varepsilon$ yields that
\begin{eqnarray}\label{3.1}
\partial_tu_\varepsilon(t,x)+b(t,x)\cdot\nabla u_\varepsilon(t,x)
+\sum_{i=1}^d\partial_{x_i}u_\varepsilon(t,x)\circ\dot{B}_i(t)=r_\varepsilon,
\end{eqnarray}
with
\begin{eqnarray*}
r_\varepsilon=b(t,x)\cdot\nabla u_\varepsilon(t,x)-(b\cdot\nabla u)_\varepsilon(t,x).
\end{eqnarray*}
With the help of assumption (\ref{1.4}), for almost all $\omega\in\Omega$,
\begin{eqnarray}\label{3.2}
\left\{
\begin{array}{ll}
r_\varepsilon\rightarrow 0  \ \ in \ \ L^q(0,T;L^p_{loc}({\mathbb R}^d)), \ \ if  \ \ p<\infty , \\ r_\varepsilon\rightarrow 0 \ \ in \ \ L^q(0,T;L^r_{loc}({\mathbb R}^d)), \ \ \forall \ 1<r<\infty,  \ if \ p=\infty.
  \end{array}
\right.
\end{eqnarray}
From (\ref{3.1}),   $u_\varepsilon(t,X(t,x))$  satisfies  that
\begin{eqnarray}\label{3.3}
\frac{d}{dt}u_\varepsilon(t,X(t,x))=r_\varepsilon(t,X(t,x)),
\end{eqnarray}
which implies that, for every $t\in (0,T)$, every $\varphi\in \mathcal{C}_0^\infty({\mathbb R}^d)$,
\begin{eqnarray}\label{3.4}
\int_{{\mathbb R}^d}u_\varepsilon(t,X(t,x))\varphi(x)dx&=&
\int^t_0\int_{{\mathbb R}^d}r_\varepsilon(s,X(s,x))\varphi(x)dx
\cr\cr&=&\int^t_0\int_{{\mathbb R}^d}r_\varepsilon(s,x))\varphi(X^{-1}(s,x))
det(\nabla_xX^{-1}(s,x))dxds.
\end{eqnarray}
In view of Lemma \ref{lem2.1}, $\varphi(X^{-1})det(\nabla_xX^{-1})$ is bounded which has a compact support. By (\ref{3.2}), if one takes $\varepsilon$ approaching to $0$, from (\ref{3.4}), one arrives at
\begin{eqnarray}\label{3.5}
\int_{{\mathbb R}^d}u(t,X(t,x))\varphi(x)dx=0,
\end{eqnarray}
which suggests that $u=0$ if one uses Lemma \ref{lem2.1} again. From this one proves the uniqueness.

Second, we show that $u(t,x)=u_0(X^{-1}(t,x))$ is a stochastic weak solution of (\ref{1.1}). Here $\{X(t,x)\}$ is the unique
strong solution of (\ref{2.1}) with $s=0$. By Lemma \ref{lem2.1}, the stochastic differential equation
(\ref{2.1}) with $s=0$ has a unique continuous adapted solution $\{X(t,x), \ t\in [0,T ], \ \omega \in\Omega\}$, which forms a $\mathcal{C}^{1,\alpha^\prime}$ ($\alpha^\prime<\alpha-2/q$) stochastic flow $X(t,x)$ of diffeomorphisms. If one defines $u(t,x)=u_0(X^{-1}(t,x))$ and uses the Kunita-It\^{o}-Wentzel formula (see \cite[Theorem 8.3]{Kun1} or \cite[Lemma 2.1]{CO}), for every $\varphi\in\mathcal{C}_0^\infty({\mathbb R}^d)$, $\int_{{\mathbb R}^d}\varphi(x)u(t,x)dx$ meets (\ref{1.2}). Thus $\int_{{\mathbb R}^d}\varphi(x)u(t,x)dx$ has a continuous modification which is an $\mathcal{F}_t$-semimartingale. To
complete the proof, we need to show $u\in
L^\infty(\Omega\times(0,T);L^p({\mathbb R}^d))$. Clearly when $p=\infty$, it is true. It remains to show $p\in [1,\infty)$.

With the help of Euler's identity,  we have
\begin{eqnarray}\label{3.6}
\int_{{\mathbb R}^d}|u_0(X^{-1}(t,x))|^pdx&=&\int_{{\mathbb R}^d}|u_0(x)|^pdet(\nabla_xX(t,x))dx
\cr\cr&=&\int_{{\mathbb R}^d}|u_0(x)|^p\exp(\int^t_0\div b(r,X(r,x))dr)dx\cr\cr&\leq&
\exp(\|\div b\|_{L^1(0,T;L^\infty({\mathbb R}^d))})\int_{{\mathbb R}^d}|u_0(x)|^pdx.
\end{eqnarray}

Third, we show that (\ref{1.3}) holds. Noticing that the stochastic differential equation (\ref{2.1}) beginning from $s=0$ has a unique continuous adapted solution $\{X(t,x), \ t\in [0,T ], \ \omega \in\Omega\}$, which forms a $\mathcal{C}^{1,\alpha^\prime}$ ($0<\alpha^\prime<\alpha-2/q$) stochastic flow $X(t,x)$ of diffeomorphisms.
We have the following chain rule
\begin{eqnarray}\label{3.7}
\nabla_x(u_0(X^{-1}(t,x)))=\nabla_xu_0(X^{-1}(t,x))
\nabla_xX^{-1}(t,x).
\end{eqnarray}
Recall the proof of \cite[Theorem 1.1 (i)]{WDGL} (or see Appendix Lemma A.1), for every $r\in [1,\infty)$, every $R>0$,
\begin{eqnarray}\label{3.8}
{\mathbb E}\sup_{0\leq t\leq T, x\in B_R}\|\nabla_xX^{-1}(t,x)\|^r\leq C(T,d,r,R)<\infty.
\end{eqnarray}
Combining (\ref{3.7}) and (\ref{3.8}), when $p<\infty$,
\begin{eqnarray}\label{3.9}
\int_{B_R}|\nabla_x(u_0(X^{-1}(t,x)))|^pdx\leq
\int_{{\mathbb R}^d}|\nabla_xu_0(x)|^pdx {\mathbb E}\sup_{0\leq t\leq T, x\in B_R}\|\nabla_xX^{-1}(t,x)\|^p<\infty,
\end{eqnarray}
and when $p=\infty$,
\begin{eqnarray}\label{3.10}
{\mathbb E}\|\nabla_x(u_0(X^{-1}))\|_{L^\infty((0,T)\times B_R)}^r\leq
\|\nabla_xu_0\|^r_{L^\infty({\mathbb R}^d)}{\mathbb E}\sup_{0\leq t\leq T, x\in B_R}\|\nabla_xX^{-1}(t,x)\|^r<\infty,
\end{eqnarray}
for every $r\in [1,\infty)$. From estimates (\ref{3.9}), (\ref{3.10}), one accomplishes the proof. $\Box$

\vskip2mm\noindent
\textbf{Remark 3.1.} From our proof one also asserts that: if $u_0\in \mathcal{C}^1_b({\mathbb R}^d)$, $b\in L^\infty(0,T;\mathcal{C}_b^\beta({\mathbb R}^d;{\mathbb R}^d))$, there exists a unique classical $\mathcal{C}^1$-solution $u$ of the Cauchy problem (\ref{1.1}), i.e.  $u(\omega,t)\in \mathcal{C}^1({\mathbb R}^d)$, for every $\varphi\in\mathcal{C}_0^\infty({\mathbb R}^d)$, $\int_{{\mathbb R}^d}\varphi(x)u(t,x)dx$ has a continuous modification which is an $\mathcal{F}_t$-semimartingale and  for every $t\in [0,T]$, (\ref{1.2}) holds. Moreover, $u(t,x)=u_0(X^{-1}(t,x))$.  We generalize the result \cite[Theorem 7]{FGP3}, which is given by Flandoli, Gubinelli and Priola.

\section{Proof of Theorem \ref{the1.2}}\label{sec4}
\setcounter{equation}{0}
\vskip2mm\noindent
\textbf{Proof.}  Now let us examine the non-existence of stochastic strong solutions. Without loss of generality, we suppose that $d=2$ and now we rewrite $x$ by $(x,y)\in{\mathbb R}^2$. We divide the proof into two cases.
 \vskip2mm\par
\textbf{$\bullet$ Case 1:} $p\in[1,\infty)$.
\vskip2mm\par
For $\epsilon$ being a small enough positive real number, we define $f(t)$ and $g(x)$ as the following:
\begin{eqnarray}\label{4.1}
g(x)=\left\{ \begin{array}{ll}
x^{\alpha},  \quad 0<x<1, \\ \ 1, \quad \ 1\leq x, \\ \ 0, \quad \ x\leq 0,\end{array}
\right.
f(t)=\left\{ \begin{array}{ll}
(t_1-t)^{-\frac{\alpha+1}{2}-\epsilon},  \quad 0\leq t<t_1\leq T, \\ \quad \quad \ 0, \quad \quad \quad\quad \ otherwise.\end{array}
\right.
\end{eqnarray}
Then $0\leq fg\in L^q(0,T;\mathcal{C}_b^\alpha({\mathbb R}))$ with $1+\alpha<2/q$. We define $b(t,x,y)=(0,f(t)g(x))$, then $b \in L^q(0,T;\mathcal{C}_b^\alpha({\mathbb R}^2))$ and
$\div b(t,x,y)=0$. Consider the SDE below
\begin{eqnarray}\label{4.2}
dX(t)=dB_1(t), \ dY(t)=f(t)g(X(t))dt+dB_2(t),
\ X(0)=x, \ Y(0)=y.
\end{eqnarray}
We obtain that
\begin{eqnarray}\label{4.3}
X(t,x)=x+B_1(t), \ Y(t,x,y)=y+B_2(t)+\int_0^tf(s)g(x+B_1(s))ds,
\end{eqnarray}
which hints that the solution of SDE (\ref{4.2}) forms a stochastic flow of homeomorphism. And from (\ref{4.3}),
\begin{eqnarray*}
\frac{\partial(X,Y)}{\partial(x,y)}=
\left(  \begin{array}{cc} 1 & \int_0^tf(s)g^\prime(x+B_1(s))ds \\                                        0 & 1 \\
\end{array}  \right),
\Big(\frac{\partial(X,Y)}{\partial(x,y)}\Big)^{-1}=
\left(  \begin{array}{cc} 1 & -\int_0^tf(s)g^\prime(x+B_1(s))ds \\                                        0 & 1 \\
\end{array}  \right).
\end{eqnarray*}
Hence the stochastic flow is measure-preserving.

For $u_0\in W^{1,p}({\mathbb R}^2)$ with $p\in [1,\infty)$, we define $u(t,x,y)=u_0((X,Y)^{-1}(t,x,y))$. Following the discussion in proving Theorem \ref{the1.1} (second step), $u(t,x)$ is a stochastic weak solution of (\ref{1.1}). Now we inspect the uniqueness.

Let $\varrho_1$ be a regularizing kernel on ${\mathbb R}$. For $\varepsilon_1>0, \varepsilon_2>0$,  denote
\begin{eqnarray*}
\varrho_{1,\varepsilon_1} =\frac{1}{\varepsilon_1} \varrho_1(\frac{\cdot}{\varepsilon_1}), \ \varrho_{1,\varepsilon_2} =\frac{1}{\varepsilon_2} \varrho_1(\frac{\cdot}{\varepsilon_2}).
\end{eqnarray*}
We define $u_{\varepsilon_2}(t,x,y)=(u(t,x,\cdot)\ast \varrho_{1,\varepsilon_2})(y)$ and $u_{\varepsilon_1,\varepsilon_2}(t,x,y)=(u\ast \varrho_{1,\varepsilon_2}(t,\cdot,y))\ast \varrho_{1,\varepsilon_1}(x)=:u_\varepsilon(t,x,y)$.  Then $u_\varepsilon$ yields that
\begin{eqnarray*}
\partial_tu_\varepsilon(t,x,y)+\partial_xu_\varepsilon(t,x,y)+
f(t)g(x)\partial_yu_\varepsilon(t,x,y)
+\nabla_{x,y}u_\varepsilon(t,x,y)\cdot\circ\dot{B}(t)=r_\varepsilon(t,x,y),
\end{eqnarray*}
with
\begin{eqnarray*}
r_\varepsilon=f(t)g(x)\partial_yu_\varepsilon(t,x,y)-f(t)
(g(\cdot)\partial_yu_{\varepsilon_2}(t,\cdot,y))\ast
\varrho_{1,\varepsilon_1}(x).
\end{eqnarray*}
Clearly, for $\varepsilon_2>0$ be fixed,  for almost all $\omega\in\Omega$,
\begin{eqnarray*}
r_\varepsilon\rightarrow 0  \ \ in \ \ L^1(0,T;L^1_{loc}({\mathbb R}^2)).
\end{eqnarray*}
Repeating the calculations from (\ref{3.3}) to (\ref{3.4}),  and    taking $\varepsilon_1$ to $0$ first, $\varepsilon_2$ to $0$ next, we arrive  at
the identity (\ref{3.5}). So the stochastic weak solution is unique.

Owing to the explicit expression of $X$ and $Y$ given by (\ref{4.3}), we have
\begin{eqnarray}\label{4.4}
X^{-1}(t,x)=x-B_1(t), \ Y^{-1}(t,x,y)=y-B_2(t)-\int_0^tf(s)g(x-B_1(s))ds.
\end{eqnarray}

Let $t_1>0$ be given in (\ref{4.1}). Then  for every $R>0$,
\begin{eqnarray}\label{4.5}
&&\int_{[-R,R]^2}|\nabla_{x,y}(u_0((X,Y)^{-1}(t_1,x,y)))|^pdxdy\cr\cr&=&
\int_{[-R,R]^2}|\nabla_{x,y}u_0((X,Y)^{-1}(t_1,x,y))|^p
\Big\|\Big(\frac{\partial(X,Y)^{-1}}{\partial(x,y)}\Big)\Big
\|^pdxdy
\cr\cr&\geq& \int_{[-R,R]^2}|\nabla_{x,y}u_0((X,Y)^{-1}(t_1,x,y))|^p
\Big|\int_0^{t_1}f(s)g^\prime(x-B_1(s))ds\Big|^pdxdy
\cr\cr&=& \int_{(X,Y)^{-1}(t_1)([-R,R]^2)}|\nabla_{x,y}u_0(x,y)|^p
\Big|\int_0^{t_1}f(s)g^\prime(x+B_1(t_1)-B_1(s))ds\Big|^pdxdy.
\end{eqnarray}

In view of (\ref{4.4}),
\begin{eqnarray}\label{4.6}
&&(X,Y)^{-1}(t_1)([-R,R]^2)\cr\cr&=&\{(x-B_1(t_1), \ y-B_2(t_1)-\int_0^{t_1}f(s)g(x-B_1(s))ds), \ (x,y)\in [-R,R]^2\}.
\end{eqnarray}

On the other hand, by (\ref{4.1}),  we observe  that
\begin{eqnarray}\label{4.7}
\int_0^{t_1}f(s)g(x-B_1(s))ds \leq \frac{2}{1-\alpha-2\epsilon}t_1^{\frac{1-\alpha-2\epsilon}{2}}\leq \frac{2}{1-\alpha-2\epsilon}T^{\frac{1-\alpha-2\epsilon}{2}}.
\end{eqnarray}
Choosing $R>6/(1-\alpha-2\epsilon)T^{\frac{1-\alpha-2\epsilon}{2}}$,  and noticing  (\ref{4.6}) and (\ref{4.7}), it follows that
\begin{eqnarray}\label{4.8}
(X,Y)^{-1}(t_1)([-R,R]^2)\supset[-R-B_1(t_1),R-B_1(t_1)]\times [-R-B_2(t_1),
\frac{2R}{3}-B_2(t_1)].
\end{eqnarray}

By (\ref{4.8}), from (\ref{4.5})  we have
\begin{eqnarray}\label{4.9}
&&\int_{[-R,R]^2}|\nabla_{x,y}(u_0((X,Y)^{-1}(t_1,x,y)))|^pdxdy
\cr\cr&\geq&
\int_{-R-B_1(t_1)}^{R-B_1(t_1)}|u_{0,1}^\prime(x)|^p
\Big|\int_0^{t_1}f(s)g^\prime(x+B_1(t_1)-B_1(s))ds\Big|^pdx
\int_{-R-B_2(t_1)}^{\frac{2R}{3}-B_2(t_1)}|u_{0,2}^\prime(y)|^pdy,
\end{eqnarray}
if one fetches $u_0(x,y)=u_{0,1}(x)u_{0,2}(y)$.

For every $s_1,s_2>0$, $B_1(s_1)$ and $B_2(s_2)$ are independent, which have the normal distributions with expected value $0$, variance $s_1$, and expected value $0$, variance $s_2$ respectively. Let the two events $\Omega_1$ and $\Omega_2$ be defined by the following
\begin{eqnarray}\label{4.10}
\Omega_1=\{\omega\in\Omega, \ |B_1(t_1)|\leq \frac{R}{3}\}, \quad \Omega_2=\{\omega\in\Omega, \ |B_2(t_1)|\leq \frac{R}{3}\}.
\end{eqnarray}
Then $\Omega_1$ and $\Omega_2$ are independent and ${\mathbb P}(\Omega_1)={\mathbb P}(\Omega_2)>0$. Therefore,
\begin{eqnarray}\label{4.11}
&&{\mathbb E}\int_{[-R,R]\times [-R,R]}|\nabla_{x,y}(u_0((X,Y)^{-1}(t_1,x,y)))|^pdxdy
\cr\cr&\geq&
{\mathbb E} \Big[ 1_{\Omega_1}1_{\Omega_2}\int_{-R-B_1(t_1)}^{R-B_1(t_1)}|u_{0,1}^\prime(x)|^p
\Big|\int_0^{t_1}f(s)g^\prime(x+B_1(t_1)-B_1(s))ds\Big|^pdx  \int_{-R-B_2(t_1)}^{\frac{2R}{3}-B_2(t_1)}|u_{0,2}^\prime(y)|^pdy\Big]
\cr\cr&=&
{\mathbb E} \Big[ 1_{\Omega_1}\int_{-R-B_1(t_1)}^{R-B_1(t_1)}|u_{0,1}^\prime(x)|^p
\Big|\int_0^{t_1}f(s)g^\prime(x+B_1(t_1)-B_1(s))ds\Big|^pdx\Big]
\cr\cr&&\times{\mathbb E} \Big[1_{\Omega_2}\int_{-R-B_2(t_1)}^{\frac{2R}{3}-B_2(t_1)}|u_{0,2}^\prime(y)|^pdy\Big]
\cr\cr&\geq&
{\mathbb E} \Big[ 1_{\Omega_1}\int_{-R-B_1(t_1)}^{R-B_1(t_1)}|u_{0,1}^\prime(x)|^p
\Big|\int_0^{t_1}f(s)g^\prime(x+B_1(t_1)-B_1(s))ds\Big|^pdx\Big]
\Big[\int_{-\frac{2R}{3}}^{\frac{R}{3}}|u_{0,2}^\prime(y)|^pdy\Big]
\cr\cr&=&C
{\mathbb E} \Big[ 1_{\Omega_1}\int_{-R-B_1(t_1)}^{R-B_1(t_1)}|u_{0,1}^\prime(x)|^p
\Big|\int_0^{t_1}f(s)g^\prime(x+B_1(t_1)-B_1(s))ds\Big|^pdx\Big]
\cr\cr&\geq&C
 \int_{-\frac{2R}{3}}^{\frac{2R}{3}}|u_{0,1}^\prime(x)|^p
\Big|\int_0^{t_1}f(s){\mathbb E} \Big[ 1_{\Omega_1}g^\prime(x+B_1(t_1)-B_1(s))\Big]ds\Big|^pdx
\cr\cr&\geq&C
 \int_{-\frac{R}{8}}^{\frac{R}{8}}|u_{0,1}^\prime(x)|^p
\Big|\int_0^{t_1}f(s){\mathbb E} \Big[ 1_{\Omega_1}g^\prime(x+B_1(t_1)-B_1(s))\Big]ds\Big|^pdx .
\end{eqnarray}
By (\ref{4.10}),  we infer that
\begin{eqnarray}\label{4.12}
&&{\mathbb E} \Big[ 1_{\Omega_1}g^\prime(x+B_1(t_1)-B_1(s))\Big]
\cr\cr&=&
{\mathbb E} \Big[ 1_{(-\frac{R}{3},\frac{R}{3})}(B_1(t_1))g^\prime(x+B_1(t_1)-B_1(s))\Big]
\cr\cr&=&
{\mathbb E} \Big[ {\mathbb E}[ 1_{(-\frac{R}{3},\frac{R}{3})}(B_1(t_1)-B_1(s)+B_1(s))g^\prime(x+B_1(t_1)-B_1(s))|\mathcal{F}_s\Big]
\cr\cr&=&
\int_{{\mathbb R}}\frac{1}{\sqrt{2\pi s}} e^{-\frac{y^2}{2s}}dy
\int_{-\frac{R}{3}}^{\frac{R}{3}}\frac{1}{\sqrt{2\pi (t_1-s)}}e^{-\frac{(z-y)^2}{2(t_1-s)}}g^\prime(x+z-y)dz
\cr\cr&\geq&Cs^{-\frac{1}{2}}(t_1-s)^{-\frac{1}{2}}e^{-\frac{x^2}{(t_1-s)}}
\int_{{\mathbb R}}e^{-\frac{y^2}{2s}}dy
\int_{-\frac{R}{3}}^{\frac{R}{3}}e^{-\frac{(z-y+x)^2}{(t_1-s)}}g^\prime(x+z-y)dz
\cr\cr&\geq&Cs^{-\frac{1}{2}}(t_1-s)^{-\frac{1}{2}}e^{-\frac{x^2}{(t_1-s)}}
\int_{-\frac{R}{8}}^{\frac{R}{8}}e^{-\frac{y^2}{2s}}dy
\int_{-\frac{R}{3}-y+x}^{\frac{R}{3}-y+x}e^{-\frac{z^2}{(t_1-s)}}g^\prime(z)dz.
\end{eqnarray}
Noting that in (\ref{4.11}) $x\in [-R/8,R/8]$, so $-y+x\in [-R/4,R/4]$ for $y\in [-R/8,R/8]$. From (\ref{4.12}) one reaches at
\begin{eqnarray}\label{4.13}
&&{\mathbb E} \Big[ 1_{\Omega_1}g^\prime(x+B_1(t_1)-B_1(s))\Big]
\cr\cr&\geq&Cs^{-\frac{1}{2}}(t_1-s)^{-\frac{1}{2}}e^{-\frac{x^2}{(t_1-s)}}
\int_{-\frac{R}{8}}^{\frac{R}{8}}e^{-\frac{y^2}{2s}}dy
\int_{-\frac{R}{12}}^{\frac{R}{12}}e^{-\frac{z^2}{(t_1-s)}}g^\prime(z)dz.
\end{eqnarray}

The function $g$ is given by (\ref{4.1}), thus
\begin{eqnarray}\label{4.14}
\int_{-\frac{R}{12}}^{\frac{R}{12}}e^{-\frac{z^2}{(t_1-s)}}g^\prime(z)dz
&=&(t_1-s)^{\frac{\alpha}{2}}
\int_{-\frac{R}{12\sqrt{t_1-s}}}^{\frac{R}{12\sqrt{t_1-s}}}g^\prime(y)e^{-y^2}dy
\cr\cr&\geq& C(t_1-s)^{\frac{\alpha}{2}}
\int_{-\frac{R}{12\sqrt{T}}}^{\frac{R}{12\sqrt{T}}}g^\prime(y)e^{-y^2}dy
\cr\cr&=& C(t_1-s)^{\frac{\alpha}{2}}.
\end{eqnarray}
Moreover,
\begin{eqnarray}\label{4.15}
\int_{-\frac{R}{8}}^{\frac{R}{8}}e^{-\frac{y^2}{2s}}dy=s^{\frac{1}{2}}
\int_{-\frac{R}{8\sqrt{s}}}^{\frac{R}{8\sqrt{s}}}e^{-\frac{y^2}{2}}dy
\geq s^{\frac{1}{2}}\int_{-\frac{R}{8\sqrt{T}}}^{\frac{R}{8\sqrt{T}}}e^{-\frac{y^2}{2}}dy
=Cs^{\frac{1}{2}}.
\end{eqnarray}

Combining (\ref{4.14}) and (\ref{4.15}), from (\ref{4.13}), we fulfill that
\begin{eqnarray}\label{4.16}
{\mathbb E} \Big[ 1_{\Omega_1}g^\prime(x+B_1(t_1)-B_1(s))\Big]
\geq C(t_1-s)^{\frac{\alpha-1}{2}}e^{-\frac{x^2}{(t_1-s)}}.
\end{eqnarray}
With the aid of (\ref{4.1}) and (\ref{4.16}),  we conclude that
\begin{eqnarray}\label{4.17}
\int_0^{t_1}f(s){\mathbb E} \Big[ 1_{\Omega_1}g^\prime(x+B_1(t_1)-B_1(s))\Big]ds&\geq &C\int_0^{t_1}f(s)(t_1-s)^{\frac{\alpha-1}{2}}e^{-\frac{x^2}{(t_1-s)}}ds
\cr\cr&=&C \int_0^{t_1}s^{-\frac{\alpha+1}{2}-\epsilon} e^{-\frac{x^2}{s}}
s^{\frac{\alpha-1}{2}}ds
\cr\cr&=&C x^{-2\epsilon}\int_{\frac{x^2}{t_1}}^\infty s^{\epsilon-1} e^{-s}ds.
\end{eqnarray}

By (\ref{4.17}), if $u_{0,1}^\prime(x)\approx x^{\epsilon-\frac{1}{p}}$ near $0+$,  and noticing  (\ref{4.11}), we obtain that
\begin{eqnarray}\label{4.18}
&&{\mathbb E}\int_{[-R,R]\times [-R,R]}|\nabla_{x,y}(u_0((X,Y)^{-1}(t_1,x,y)))|^pdxdy
\cr\cr&\geq&C
\int_0^{\frac{R}{8}}x^{p\epsilon-1}
\Big|x^{-2\epsilon}\int_{\frac{x^2}{t_1}}^\infty s^{\epsilon-1} e^{-s}ds\Big|^pdx
\cr\cr&\geq& C
\int_0^{\frac{R}{8}}x^{-p\epsilon-1}\Big|\int_{\frac{R^2}{64T}}^\infty s^{\epsilon-1} e^{-s}ds\Big|^pdx
\cr\cr&\geq & C\int_0^{\frac{R}{8}}x^{-p\epsilon-1} dx=\infty.
\end{eqnarray}

\vskip2mm\par
\textbf{$\bullet$ Case 2:} $p=\infty$.
\vskip2mm\par
Let $\epsilon$ be a small enough positive real number. We define $g(x)$ by (\ref{4.1}), $h(t)$:
\begin{eqnarray}\label{4.19}
h(t)=\left\{ \begin{array}{ll}
t^{-\frac{\alpha+1}{2}-\epsilon},  \quad 0\leq t\leq T, \\ \quad \ \ 0, \quad \quad \ otherwise.\end{array}
\right.
\end{eqnarray}
Then $0\leq hg\in L^q(0,T;\mathcal{C}_b^\alpha({\mathbb R}))$ with $1+\alpha<2/q$. Consider the SDE below
\begin{eqnarray*}
dX(t)=dB_1(t), \ dY(t)=h(t)g(X(t))dt+dB_2(t).
\ X(0)=x, \ Y(0)=y,
\end{eqnarray*}
Hence
\begin{eqnarray*}
X(t,x)=x+B_1(t), \ Y(t,x,y)=y+B_2(t)+\int_0^th(s)g(x+B_1(s))ds.
\end{eqnarray*}
For $u_0\in W^{1,\infty}({\mathbb R}^2)$, we define $u(t,x,y)=u_0((X,Y)^{-1}(t,x,y))$. Following the discussion in Case 1, $u(t,x)$ is the unique stochastic weak solution of (\ref{1.1}). Now we show the regularity.

In view of (\ref{4.4}), (\ref{4.7}) and (\ref{4.8}), for every $R>6/(1-\alpha-2\epsilon)T^{\frac{1-\alpha-2\epsilon}{2}}>0$,  and taking
$u_0(x,y)=u_{0,1}(x)u_{0,2}(y)$, for every $r\geq 1$ and every $t\in (0,T)$, we have
\begin{eqnarray}\label{4.20}
&&{\mathbb E}\|\nabla_{x,y}(u_0((X,Y)^{-1}(t)))\|^r_{L^\infty([-R,R]^2)}
\cr\cr&\geq&
{\mathbb E}\Big[1_{\Omega_3} \sup_{x\in (-R,R)}\Big|u_{0,1}^\prime(x-B_1(t))\int_0^{t}h(s)g^\prime(x-B_1(s))ds\Big|^r\Big]
\cr\cr&&\times
{\mathbb E}\Big[1_{\Omega_4} \sup_{y\in (-R,\frac{2R}{3})}|u_{0,2}^\prime(y-B_2(t))|^r\Big]
\cr\cr&\geq&
C\|u_{0,2}^\prime\|^r_{L^\infty(-\frac{2R}{3},\frac{R}{3})}
{\mathbb E}\Big[1_{\Omega_3} \sup_{x\in (-R,R)}\Big|u_{0,1}^\prime(x-B_1(t))\int_0^{t}h(s)g^\prime(x-B_1(s))ds\Big|^r\Big]
\cr\cr&\geq&C{\mathbb E}\Big[1_{\Omega_3} \sup_{x\in (-R,R)}\Big|u_{0,1}^\prime(x-B_1(t))\int_0^{t}h(s)g^\prime(x-B_1(s))ds\Big|^r\Big],
\end{eqnarray}
where the events $\Omega_3$ and $\Omega_4$ are defined by
\begin{eqnarray*}
\Omega_3=\{\omega\in\Omega, \ |B_1(t)|\leq \frac{R}{3}\}, \quad \Omega_4=\{\omega\in\Omega, \ |B_2(t)|\leq \frac{R}{3}\}.
\end{eqnarray*}

Let $u_{0,1}\in W^{1,\infty}({\mathbb R})$ be such that $u_{0,1}^\prime(x)=1$ when $x\in [-2R,2R]$. Then by (\ref{4.20}) we gain that
\begin{eqnarray}\label{4.21}
&&{\mathbb E}\|\nabla_{x,y}(u_0((X,Y)^{-1}(t)))\|^r_{L^\infty([-R,R]^2)}
\cr\cr&\geq&C{\mathbb E}\Big[1_{\Omega_1} \sup_{x\in (-R,R)}\Big|\int_0^{t}h(s)g^\prime(x-B_1(s))ds\Big|^r\Big]
\cr\cr&\geq&C\sup_{x\in (-R,R)}\Big|\int_0^{t}h(s){\mathbb E}[1_{\Omega_1}g^\prime(x-B_1(s))]ds\Big|^r.
\end{eqnarray}

Observe that
\begin{eqnarray}\label{4.22}
&&{\mathbb E} \Big[ 1_{\Omega_1}g^\prime(x-B_1(s))\Big]
\cr\cr&=&
{\mathbb E} \Big[ 1_{(-\frac{R}{3},\frac{R}{3})}(B_1(t))g^\prime(x-B_1(s))\Big]
\cr\cr&=&
{\mathbb E} \Big[ {\mathbb E}[ 1_{(-\frac{R}{3},\frac{R}{3})}(B_1(t)-B_1(s)+B_1(s))g^\prime(x-B_1(s))|\mathcal{F}_s\Big]
\cr\cr&=&
\int_{{\mathbb R}}\frac{1}{\sqrt{2\pi s}}g^\prime(x-y) e^{-\frac{y^2}{2s}}dy
\int_{-\frac{R}{3}}^{\frac{R}{3}}\frac{1}{\sqrt{2\pi (t-s)}}e^{-\frac{(z-y)^2}{2(t-s)}}dz
\cr\cr&\geq&C
\int_{{\mathbb R}}\frac{1}{\sqrt{2\pi s}}g^\prime(x-y) e^{-\frac{y^2}{2s}}e^{-\frac{y^2}{(t-s)}} dy
\cr\cr&\geq&Ce^{-\frac{x^2}{s}}e^{-\frac{2x^2}{(t-s)}}
\int_{{\mathbb R}}\frac{1}{\sqrt{2\pi s}}g^\prime(y) e^{-\frac{y^2}{s}}e^{-\frac{2y^2}{(t-s)}} dy.
\end{eqnarray}
Noting that the function $g$ is (\ref{4.1}),  we thus have
\begin{eqnarray}\label{4.23}
\int_{{\mathbb R}}\frac{1}{\sqrt{2\pi s}}g^\prime(y) e^{-\frac{y^2}{s}}e^{-\frac{2y^2}{(t-s)}} dy
&=&\Big(\frac{s(t-s)}{(t+s)}\Big)^{\frac{\alpha-1}{2}}
\Big(\frac{(t-s)}{(t+s)}\Big)^{\frac{1}{2}}
\int_{{\mathbb R}}g^\prime(y)e^{-y^2}dy
\cr\cr&\geq &Cs^{\frac{\alpha-1}{2}}(t-s)^{\frac{\alpha}{2}}t^{-\frac{\alpha}{2}}.
\end{eqnarray}
Combining (\ref{4.22}) and (\ref{4.23}),  we reach at
\begin{eqnarray}\label{4.24}
\int_0^th(s){\mathbb E} \Big[ 1_{\Omega_1}g^\prime(x-B_1(s))\Big]ds
&\geq &C\int_0^th(s)e^{-\frac{x^2}{s}}e^{-\frac{2x^2}{(t-s)}}
s^{\frac{\alpha-1}{2}}(t-s)^{\frac{\alpha}{2}}t^{-\frac{\alpha}{2}}ds
\cr\cr&\geq &Ce^{-\frac{4x^2}{t}} \int_0^{\frac{t}{2}}h(s)e^{-\frac{x^2}{s}}
s^{\frac{\alpha-1}{2}}ds.
\end{eqnarray}

By virtue of (\ref{4.19}),
\begin{eqnarray}\label{4.25}
\int_0^{\frac{t}{2}}h(s) e^{-\frac{x^2}{s}}s^{\frac{\alpha-1}{2}}ds= \int_0^{\frac{t}{2}}s^{-\frac{\alpha+1}{2}-\epsilon} e^{-\frac{x^2}{s}}s^{\frac{\alpha-1}{2}}ds= x^{-2\epsilon}\int_{\frac{2x^2}{t}}^\infty s^{\epsilon-1} e^{-s}ds.
\end{eqnarray}
Therefore
\begin{eqnarray*}
{\mathbb E}\|\nabla_{x,y}(u_0((X,Y)^{-1}(t)))\|^r_{L^\infty([-R,R]^2)}
\geq C \sup_{x\in (-R,R)}\Big|x^{-2\epsilon}\int_{\frac{2x^2}{t}}^\infty s^{\epsilon-1} e^{-s}ds|^r=\infty,
\end{eqnarray*}
which implies that
\begin{eqnarray*}
{\mathbb E}\|\nabla_{x,y}(u_0((X,Y)^{-1}))\|^r_{L^\infty((0,T)\times [-R,R]^2)}
\geq \sup_{0\leq t\leq T}{\mathbb E}\|\nabla_{x,y}(u_0((X,Y)^{-1}(t)))\|^r_{L^\infty([-R,R]^2)}
=\infty.
\end{eqnarray*}
So we finish the proof. $\Box$

\vskip2mm\noindent
\textbf{Remark 4.1.} When the  `Ladyzhenskaya-Prodi-Serrin type condition' (\ref{1.8}) does not hold,  and use  Theorem \ref{the1.2},  we  assert   the non-existence of stochastic strong solutions. On the other hand, from Remark 1.2,  we may find  the existence and uniqueness of $\cap_{r\geq 1}W^{1,r}_{loc}({\mathbb R}^d)$ solutions under  the `Ladyzhenskaya-Prodi-Serrin type condition'. However,  we do not know how to establish the existence and uniqueness of strong solutions for (\ref{1.1}) under the  `Ladyzhenskaya-Prodi-Serrin type condition' (\ref{1.1}). Perhaps, there is a real number $\alpha\leq \alpha_0\leq 1+\alpha$, such that for $2/q\leq \alpha_0$, Theorem \ref{the1.1} is true and when $2/q> \alpha_0$, Theorem \ref{the1.2} is true.


\begin{appendix}\label{appendixA}
\section{Appendix:  A useful lemma}\label{}\setcounter{equation}{0}
\begin{lemma} \label{A.1}  Assume that $q, \alpha$ and $b$ be stated in Lemma \ref{lem2.1}. Suppose that $\{X(t,x), \
t\in [0,T ], \ \omega \in\Omega\}$ is the stochastic flow generated by (\ref{2.1}) with $s=0$. Then for every $R>0$, every $r\in [1,\infty)$, there is a positive constant $C(r,R)$ such that
\begin{eqnarray*}
{\mathbb E}\sup_{0\leq t\leq T, x\in B_R}\|\nabla_xX^{-1}(t,x)\|^r\leq C(T,d,r,R)<\infty.
\end{eqnarray*}
\end{lemma}
\textbf{Proof.} Since the backward flow satisfies the same SDE of the forward flow with a drift coefficient of opposite sign, to calculate ${\mathbb E}\sup_{0\leq t\leq T, x\in
B_R}\|\nabla_xX^{-1}(t,x)\|^r$ it is sufficient to estimate
${\mathbb E}\sup_{0\leq t\leq T, x\in B_R}\|\nabla_xX(t,x)\|^r$. Recall that
SDE (\ref{2.1}) is equivalent to (\ref{2.14}) and $X(t)=\gamma^{-1}(t)\circ Y(t)$.  With the help of Lemma \ref{lem2.1} (iii), it suffices to manipulate ${\mathbb E}\sup_{0\leq t\leq T, y\in B_R}\|\nabla_yY(t,y)\|^r$. By scaling and shift transformations, we only need to show ${\mathbb E}\sup_{0\leq t\leq T, x\in [0,1]^d}\|\nabla Y(t,x)\|^r$. The calculations can be divided into three steps.
\vskip2mm\par
\textbf{Step 1.} Space H\"{o}lder estimates for $Y(t)$.
\vskip2mm\par

Let $Y(t,x)$ and $Y(t,y)$ be the unique strong solution of
(\ref{2.14}) with initial data $x$ and $y$ respectively.  Define
$Y_t(x,y)=Y(t,x)-Y(t,y)$ and $\tilde{b}(t,y)=\lambda
U(t,\gamma^{-1}(t,y))$, $\sigma(t,y)=I+\nabla
U(t,\gamma^{-1}(t,y))$. Then
\begin{eqnarray*}
\left\{\begin{array}{ll} d
Y_t(x,y)=[\tilde{b}(t,Y(t,x))-\tilde{b}(t,Y(t,y))]dt
+[\sigma(t,Y(t,x))-\sigma(t,Y(t,x))]d B(t), \ t\in(0,T), \\
Y_t(x,y)|_{t=0}=x-y.
\end{array}
\right.
\end{eqnarray*}
Using the It\^{o} formula and Lemma \ref{lem2.1} (ii), for $m\geq 2$, we have
\begin{eqnarray}\label{A.1}
|Y_t(x,y)|^m&\leq&|x-y|^m+C(m)\int_0^{t}|Y_s(x,y)|^m
ds\cr\cr&&+m\int_0^t|Y_s(x,y)|^{m-2}\langle Y_s(x,y),
[\sigma(s,Y(s,x))-\sigma(s,Y(s,x))]d B(s)\rangle.
\end{eqnarray}
Therefore
\begin{eqnarray*}
{\mathbb E}|Y_t(x,y)|^m\leq|x-y|^m+C(m)\int_0^t{\mathbb E}|Y_s(x,y)|^mds,
\end{eqnarray*}
which suggests that
\begin{eqnarray}\label{A.2}
\sup_{0\leq t\leq T}{\mathbb E}|Y_t(x,y)|^m\leq C(m,T)|x-y|^m.
\end{eqnarray}

On the other hand, by virtue of the BDG inequality, from
(\ref{A.1}),  we  conclude that
\begin{eqnarray*}
{\mathbb E}\sup_{0\leq s\leq
t}|Y_s(x,y)|^m&\leq&|x-y|^m+C(m)\int_0^{t}{\mathbb E}\sup_{0\leq r\leq
s}|Y_r(x,y)|^m dr\cr\cr&&+C(m){\mathbb E}\[\int_0^t|Y_s(x,y)|^{2m}ds\]^{\frac{1}{2}}.
\end{eqnarray*}
Since (\ref{A.2}) holds for every $m\geq 2$,  we obtain by the Minkowski and Gr\"{o}nwall  inequalities,
\begin{eqnarray}\label{A.3}
{\mathbb E}\sup_{0\leq t\leq T}|Y_t(x,y)|^m\leq C(m,T)|x-y|^m.
\end{eqnarray}
From this, we  also gain
\begin{eqnarray}\label{A.4}
{\mathbb E}\sup_{0\leq t\leq T}|X_t(x,y)|^m\leq C(m,T)|x-y|^m.
\end{eqnarray}

\vskip2mm\par
\textbf{Step 2.}  H\"{o}lder estimate for $\|\nabla_xY(t,x)\|$.
\vskip2mm\par

Set $\nabla_xY(t,x)$ by $\xi_t(x)$.  Then $\xi_t(x)$ yields that
\begin{eqnarray*}
d \xi_t(x)=\lambda \nabla U(t,X(t,x))\nabla\gamma^{-1}(t,Y(t))\xi_t(x)dt+\nabla^2
U(t,X(t,x))\nabla\gamma^{-1}(t,Y(t))\xi_t(x)d B(t),
\end{eqnarray*}
and $\xi_t(x)|_{t=0}=I$.

Similar calculations from (\ref{A.1}) to (\ref{A.3}) imply that, for every $m\geq 2$,
\begin{eqnarray}\label{A.5}
\sup_{x\in{\mathbb R}^d}{\mathbb E}\sup_{0\leq t\leq T}\|\xi_t(x)\|^m\leq C.
\end{eqnarray}

If one set $\xi_t(x,y)=\xi_t(x)-\xi_t(y)$, by analogue manipulations from (\ref{A.1}) to (\ref{A.2}), it yields that
\begin{eqnarray}\label{A.6}
&&{\mathbb E}\|\xi_t(x,y)\|^m \cr\cr&\leq&C(m)\int_0^t{\mathbb E}\|\xi_s(x,y)\|^mds\cr\cr&&+
C(m){\mathbb E}\int_0^t\|\xi_s(x)\|^m[|X(s,x)-X(s,y)|^m+|X(s,x)-X(s,y)|^{(\alpha-2/q)
m}]ds,
\end{eqnarray}
as  $U\in L^\infty(0,T;\mathcal{C}^{2,\alpha-2/q}_b({\mathbb R}^d;{\mathbb R}^d))$.

With the aid of (\ref{A.4}), (\ref{A.5}),  we have from (\ref{A.6})
\begin{eqnarray*}
{\mathbb E}\|\xi_t(x,y)\|^m &\leq&C(m)\int_0^t{\mathbb E}\|\xi_s(x,y)\|^mds \cr\cr&&+
C(m)\int_0^t\Big({\mathbb E}\|\xi_s(x)\|^{2m}\Big)^{\frac{1}{2}}\Big(
{\mathbb E}[|X(s,x)-X(s,y)|^{2m}\Big)^{\frac{1}{2}}ds \cr\cr&&+
C(m)\int_0^t\Big({\mathbb E}\|\xi_s(x)\|^{2m}\Big)^{\frac{1}{2}}\Big(
{\mathbb E}|X(s,x)-X(s,y)|^{2(\alpha-2/q) m}\Big)^{\frac{1}{2}}ds
\cr\cr&\leq& C(m)\int_0^t{\mathbb E}\|\xi_s(x,y)\|^mds+
C(m,T)[|x-y|^m+|x-y|^{(\alpha-2/q) m}].
\end{eqnarray*}
Thus
\begin{eqnarray*}
\sup_{0\leq t\leq T}{\mathbb E}\|\xi_t(x,y)\|^m \leq
C(m,T)[|x-y|^m+|x-y|^{(\alpha-2/q) m}].
\end{eqnarray*}

Similar manipulations of (\ref{A.2})-(\ref{A.3}) apply again,
we  end up with
\begin{eqnarray}\label{A.7}
&&{\mathbb E}\sup_{0\leq t\leq T}\|\xi_t(x,y)\|^m\cr\cr&\leq&
C(m,T)[|x-y|^m+|x-y|^{(\alpha-2/q) m}]\leq C(m,T)|x-y|^{(\alpha-2/q) m}, \ \forall \ x,y\in [0,1]^d.
\end{eqnarray}

\vskip2mm\par
\textbf{Step 3.}  ${\mathbb E}\sup_{0\leq t\leq T, x\in
[0,1]^d}\|\nabla_xY(t,x)\|^r<\infty.$
\vskip2mm\par

To this end, we introduce a sequence of sets: $\cS_n=\{z\in
\mZ^d\ | \ z2^{-n}\in [0,1]^d\}, \ n\in \mN$. For an arbitrary
$e=(e^1,_{\cdots},e^d)\in \mZ^d$ such that $\|e\|_\infty=\max_{1\leq
i\leq d}|e^i|=1$, and every $z,z+e\in\cS_n$, we define
$\xi_z^{n,e}(t)=|\xi_t((z+e)2^{-n})-\xi_t(z2^{-n})|$. Then by
(\ref{A.7}), for every $m\geq 2$,
\begin{eqnarray*}
{\mathbb E}\sup_{0\leq t\leq T}|\xi_z^{n,e}(t)|^m\leq C(m,T)2^{-n (\alpha-2/q) m}.
\end{eqnarray*}
For any $\tau>0$ and $K>0$, one sets a number of events
$\cA_{z,\tau}^{n,e}=\{\omega\in\Omega \ | \ \sup_{0\leq t\leq
T}\xi_z^{n,e}(t)\geq K\tau^n\}$ ($z,z+e\in\cS_n$), it yields that
\begin{eqnarray*}
{\mathbb P}(\cA_{z,\tau}^{n,e})\leq \frac{{\mathbb E}\sup_{0\leq t\leq
T}|\xi_z^{n,e}(t)|^m}{K^m\tau^{mn}}\leq \frac{C(m,T)2^{-n(\alpha-2/q)
q}}{K^m\tau^{mn}}.
\end{eqnarray*}

Observing that for each $n$, the total number of the events
$\cA_{z,\tau}^{n,e}$  ($z,z+e\in\cS_n$) is not greater than
$2^{c(d)n}$. Hence the probability of the union
$\cA_\tau^n=\cup_{z,z+e\in
S_n}(\cup_{\|e\|_\infty=1}\cA_{z,\tau}^{n,e})$ fulfils  the estimate
\begin{eqnarray*}
{\mathbb P}(\cA_\tau^n)\leq C(m,T)
\frac{2^{-nm(\alpha-2/q)}}{K^m\tau^{mn}}2^{c(d)n}\leq C(m,T)K^{-m}
\Big(\frac{2^{c(d)}}{(2^(\alpha-2/q)\tau)^m}\Big)^n.
\end{eqnarray*}

  Take $\tau=2^{-(\alpha-2/q)/2}$, $m>3c(d)/(\alpha-2/q)\vee 1$. Then the
probability of the event $\cA=\cup_{n\geq1}\cA_\tau^n$ can  be
estimated  by
\begin{eqnarray}\label{A.8}
{\mathbb P}(\cA)\leq C(T,d)K^{-m}.
\end{eqnarray}

For every point $x\in [0,1]^d$, we have $x=\sum_{i=0}^\infty
e_i2^{-i}$ ($\|e_i\|_\infty\leq 1$). Denote
$x_k=\sum_{i=0}^ke_i2^{-i}$. For any $\omega\overline{\in}\cA$,
we have $|\xi_t(x_{k+1})-\xi_t(x_k)|< K\tau^{k+1}$, which suggests
that
\begin{eqnarray}\label{A.9}
|\xi_t(x)-\xi_t(x_0)|\leq \sum_{k=0}^\infty|\xi_t(x_{k+1})-\xi_t(x_k)|<
K\sum_{k=0}^\infty\tau^{k+1}\leq CK,
\end{eqnarray}
where we have fetched $\tau=2^{-(\alpha-2/q)/2}$.

Set $\xi_1=\sup_{(t,x)\in [0,T]\times [0,1]^d}|\xi_t(x)-\xi_t(x_0)|$. Then
for any $0<r<m$,
\begin{eqnarray}\label{A.10}
{\mathbb E} |\xi_1|^r=r\int_0^\infty  \lambda^{r-1}{\mathbb P}(\xi_1\geq \lambda)d\lambda
=r\int_0^{CK}\lambda^{r-1}{\mathbb P}(\xi_1\geq \lambda)d\lambda+r\int_{CK}^\infty
\lambda^{r-1}{\mathbb P}(v\geq \lambda)dr.
\end{eqnarray}
According to (\ref{A.9}), (\ref{A.8}),  and (\ref{A.10}),  we have
\begin{eqnarray*}
{\mathbb E} |\xi_1|^r\leq (CK)^r+C(T,d) r\int_{CK}^\infty
\lambda^{r-1-m}d\lambda \leq (CK)^r +C(T,d)r K^{r-m},
\end{eqnarray*}
which hints that
\begin{eqnarray*}
{\mathbb E} \sup_{(t,x)\in [0,T]\times [0,1]^d}|\xi_t(x)|^r \leq C(r)\Big[{\mathbb E} |\xi_1|^r +{\mathbb E} \sup_{0\leq t\leq T}|\xi_t(x_0)|^r\Big]\leq
C(T,d,r).
\end{eqnarray*}
This completes the proof. $\Box$
\end{appendix}

\end{document}